\documentclass[11pt]{article}
\usepackage{graphicx}
\usepackage{amsmath}
\usepackage{indentfirst}
\usepackage{textcomp}
\usepackage{amsmath,amssymb}
\usepackage{fancyhdr}
\usepackage{CJK}
\usepackage[utf8]{inputenc}
\usepackage{indentfirst}
\usepackage{geometry}
\usepackage{fancyhdr}
\usepackage{indentfirst}
\numberwithin{equation}{section}

\title{The Properties of Average Gradient in Local Region}
\author{Zhongkui Ma}
\date{}

\begin{document}

\maketitle
\thispagestyle{empty}

\begin{center}
{\bf Abstract}
\end{center}

This paper studies the average gradient over the local region of a function and constructs the homogenization function of a function.
It is found that there are some good properties about the local extreme points and the global extreme points of the function.
By using the gradient algorithm, it is more effective to use the homogenization function to find the extreme values of the function.
This method implies a method of sifting out the local extreme points of a function.


\clearpage

\pagenumbering{arabic}
\setcounter{page}{1}

\section{Introduction}

Since the gradient algorithm was discovered, the derivative of the function has become an effective and versatile tool for finding certain functional indicators.
It provides solutions for solving many complex problems, especially practical problems since the development of modern science.
At the same time, the gradient algorithm has inspired researchers to propose more innovative mathematical models, such as artificial neural networks, and even developed a discipline that relies entirely on gradient algorithms, which is surprising and anti-traditionally mathematical research.
But it does not mean that the gradient algorithm solves all problems.
Non-convex problems and local extreme points have always been obstacles to the gradient algorithm.
Many researchers have proposed corresponding solutions, especially a series of simple and effective adaptive gradient descent algorithms\cite{Rum} \cite{Pol} \cite{Duc} \cite{Die}.

The author of this paper studied the general shape of a function over a local area from the geometric point of view, i.e., the nature of the derivative or the partial derivative.
The concept of average gradient, average gradient function, homogenization function, etc., are further defined, and the case of the unary function is generalized to the multivariate function. This method is able to filter out the global extreme point(s) in the local extreme points of a function.

\section{Average Gradient and Average Gradient Function}

Consider the example that if a function is in an interval, only a small part of it decreases and the rest increases.
So it is still increasing in general.
What indicators should be used to measure this situation at this time?
This is the reason for defining the average gradient, and through further analysis, it is found that it does have a concise form and stronger nature.
First, the definition of average gradient is given.

\newtheorem{definition}{\hspace{2em}Definition}[section]
\begin{definition}
If the function $f(x)$ is continuous and derivable on the interval $[A, B]$, the average of all its derivatives on the interval is called the {\bf average derivative} of the function on the interval $[A, B]$.
For narrative convenience, it is also called the {\bf average gradient} at point $\frac{B-A}{2}$.
\end{definition}

The average gradient is relative as a measurement of the function, which means to what extent it is incremental for the examples given above.
That is the parameter involved in the definition of the average gradient and interval length.
So the following definition is made.
\newtheorem{definition3}{\hspace{2em}Definition}[section]
\begin{definition}
The value of $\left| A-B \right|$ is the {\bf scale} of the average gradient of the point $\frac{B-A}{2}$.
\end{definition}

Then, the meaning of the average gradient of the function at a certain point is defined, which describes the properties of the function in the region around this point.
\newtheorem{theorem}{\hspace{2em}Theorem}[section]
\begin{theorem}
The average gradient of the function $y=f(x)$ at point $x_0$ with scale of $h$ is $T$, and its value equals
\begin{equation}
\frac{f(x_0+\frac{h}{2}) - f(x_0-\frac{h}{2})}{h}.
\end{equation}
\end{theorem}

\newtheorem{proof}{\hspace{2em}Proof}[section]
\begin{proof}
In the interval of $[x_0-\frac{h}{2},x_0+\frac{h}{2}]$, take a column of $\{x_i\}$, then the average value of all the corresponding derivatives of these points is
$\frac{\sum_i^n f'(x_i)}{n}.$
When $n \to \infty$,
$${
\begin{aligned}
\lim\limits_{n \to \infty} \frac{\sum_i^n f'(x_i)}{n}
&= \lim\limits_{n \to \infty} \frac{\sum_i^n f'(x_i)}{n} \frac{h}{h}\\
&= \lim\limits_{n \to \infty} \sum_i^n f'(x_i) \frac{h}{n} \frac{1}{h}\\
&= \lim\limits_{n \to \infty} \sum_i^n f'(x_i) \triangle x \frac{1}{h}\\
&= \frac{1}{h} \int_{x_0-\frac{h}{2}}^{x_0+\frac{h}{2}}f'(x)dx\\
&= \frac{f(x_0+\frac{h}{2})-f(x_0-\frac{h}{2})}{h}.
\end{aligned}
}$$
\end{proof}[QED]

With the above theorem, an average gradient function corresponding to a function can be defined in the define domain of function.
At the same time, some theorems about monotone function are obtained.
\newtheorem{definition4}{\hspace{2em}Definition}[section]
\begin{definition}
For any $x_0$, the continuous and derivable function has an average gradient of $T_{x_0}$. All the $T_{x_0}$ forms a continuous derivative function, which becomes the average gradient function of the function $f(x)$ with the scale of $h$. This new function,
\begin{equation}
T(h,x)=\frac{f(x+\frac{h}{2})-f(x-\frac{h}{2})}{h},
\end{equation}
\end{definition}
is the {\bf average gradient function} of the function $f(x)$.

\newtheorem{theorem3}{\hspace{2em}Theorem}[section]
\begin{theorem}
If the average gradient function of the function $f(x)$ is positive (or negative) at any scale of $h$, then the function $f(x)$ increases (or decreases) monotonously.
\end{theorem}

\newtheorem{proof3}{\hspace{2em}Proof}[section]
\begin{proof}
Let's prove the first case of $f(x)$ increases monotonously first.

$\because$ $\forall x, h, T(h,x)=\frac{f(x+\frac{h}{2})-f(x-\frac{h}{2})}{h}>0 \text{ and } h>0.$

$\therefore$ $f(x+\frac{h}{2})-f(x-\frac{h}{2})>0.$

$\because$ $x+\frac{h}{2}>x-\frac{h}{2}.$

$\therefore$ $\frac{f(x+\frac{h}{2})-f(x-\frac{h}{2})}{x+\frac{h}{2}-(x-\frac{h}{2})}>0.$

$\therefore$ Function $f(x)$ increases monotonously.

Similarly, the case of decreasing will be proved.
\end{proof}[QED]

The proof of the above theorem has shown that the so-called average gradient function is the same as the derivative in discrete form.
This shows that the discrete derivative is not simply the estimation of the derivative. In fact, it has some interesting properties.
This is what will be demonstrated next.
This form, which has been used for a long time, has some undetected properties.

\newtheorem{theorem4}{\hspace{2em}Theorem}[section]
\begin{theorem}
If the function $f(x)$ has only two local minimum points in the interval of $[A,B]$, $x_1$ and $x_2 (x_1<x_2)$, and $f(x_1)<f(x_2)$.
There is an integer $h_0 (0<h_0<\infty)$ such that the average gradient function of $f(x)$ is $T(h, x)$ with only one zero point $x_T$ in the interval.
\end{theorem}

\newtheorem{proof4}{\hspace{2em}Proof}
\begin{proof}
$\because$ The continuous derivative $f(x)$ has two minimum points, $x_1$and $x_2$.

$\therefore$ $f'(x_1)=f'(x_2)=0,$ and $f'(x)$ is not constant $0$.

$\therefore$ In the interval $(x_1,x_2)$, there must be a point $x_3$ that is the local maximum point, $f'(x_3)=0,$ which can be obtained by the mean value theorem.

For intervals $[A,x_1]$ and $[x_2,B]$, the maximum values are denoted as $f(x_A)$ and $f(x_B)$, respectively.

The discussion is divided into the following situations:

The first situation, when $f(x_3)<f(x_A)$ and $f(x_3)<f(x_B)$,

$b \in (x_2,B], f(b)=f(x_3)$.

$\because f(x_2)>f(x_1)$.

$\therefore$ There must have $a \in (x_1,x_3)$, $f(a)=f(x_2)$.

Take $h_0=\left| a-b \right|$, then $\forall x \in [x_1+\frac{h_0}{2},x_B-\frac{h_0}{2}]$, $T(x,h_0) \ge 0$.

There must have one point $d$ in the interval $(x_2,B]$ and $f(d)=f(B)$.

Take $h_0=\left| d-B \right|$, then $\forall x \in [x_2+\frac{h_0}{2},B-\frac{h_0}{2}]$, $T(x,h_0) \ge 0$.

There must have one point $c$ in the interval $[A,x_1)$ and $f(c)=f(A)$.

Take $h_0=\left| c-A \right|$, then $\forall x \in [A+\frac{h_0}{2},x_1-\frac{h_0}{2}]$, $T(x,h_0) \le 0$.

Take $h_0=max\{\left| a-b\right|,\left| d-B\right|,\left| c-A\right|\}$.

Now, $[x_1+\frac{h_0}{2},B-\frac{h_0}{2}]=[x_1+\frac{h_0}{2},x_B-\frac{h_0}{2}] \cup [x_2+\frac{h_0}{2},B-\frac{h_0}{2}]$.

And when $x \in [x_1+\frac{h_0}{2},B-\frac{h_0}{2}]$, $T(h_0,x) \ge 0$. When $x \in [A+\frac{h_0}{2},x_1-\frac{h_0}{2}]$, $T(x,h_0) \le 0$.

$\because T(h_0,x)$ is continuous and derivable function.

$\therefore T(h_0,x)$ has a zero point in $[x_1-\frac{h_0}{2},x_1+\frac{h_0}{2}]$, and can only have 1 zero point.

If there are 2 zero points, set to $x_{T_1}$ and $x_{T_2}$, then $T(h_0,x_{T_1})=T(h_0,x_{T_2})=0$, i.e. $f (x_{T_1}+\frac{h_0}{2})-f(x_{T_1}-\frac{h_0}{2})=0$, $f(x_{T_2}+\frac{h_0}{ 2})-f(x_{T_2}-\frac{h_0}{2})=0$.

$\therefore$ $f(x_{T_1}+\frac{h_0}{2})=f(x_{T_1}-\frac{h_0}{2})$, $f(x_{T_2}+\frac {h_0}{2})=f(x_{T_2}-\frac{h_0}{2})$,
and $\left| x_{T_1}+\frac{h_0}{2}-(x_{T_1}+\frac{h_0}{2})\right| = \left| x_{T_1}+\frac{h_0 }{2}-(x_{T_1}+\frac{h_0}{2})\right|=h_0$.

$\because$ There is only one point x $\in [x_A,x_3]$.

And it satisfies $f(x+\frac{h_0}{2})=f(x-\frac{h_0}{2})$, and $\left| x+ \frac{h_0}{2}-(x-
\frac{h_0}{2}) \right|=h_0$.

$\therefore x_{T_1}$ and $x_{T_2}$ are the same point, which is $x_{T_0}$ in the theorem.

The second case, when $f(x_3)>f(x_A)$ and $f(x_3)>f(x_B)$,

Let the maximum value in the interval $[A,x_1]$ be $f(x_A)$, and the maximum value in $(x_2,B]$ is $f(x_B)$.

$\forall x \in (x_1,x_3]$, must have one point $a$, and $f(a)=f(x_A)$.

$\forall x \in [x_3,x_2)$, must have one point $b$, and $f(b)=f(x_B)$.

Take $h_0=\left| a-x_A \right|$, then $\forall x \in [x_A+\frac{h_0}{2},x_3-\frac{h_0}{2}]$, $T(h_0 ,x) \ge 0$.

Take $h_0=\left| x_B-b \right|$, then $\forall x \in [x_3+\frac{h_0}{2},x_B-\frac{h_0}{2}]$, $T(h_0 ,x) \le 0$.

Take $h_0=\left| dB \right|$, then $\forall x \in [x_2+\frac{h_0}{2},B-\frac{h_0}{2}]$, $T(h_0,x ) \le 0$.

Take $h_0=\left| cA \right|$, then $\forall x \in [A+\frac{h_0}{2},x_1-\frac{h_0}{2}]$, $T(h_0,x ) \ge 0$.

Other arguments are similar to the first case.

Other cases of $f(x_3)$, $f(x_A)$ and $f(x_B)$ are similarly proved.

\end{proof}[QED]

The proof of this theorem is the most critical. It uses the property of average gradient to screen all extreme points.
Obviously, the average gradient function corresponding to a function with many extremum points (one of which is the global maximum (or minimum) point) has only one zero point.
In order to illustrate its universality, the following corollaries are made, and their proofs are simple and repetitive.
Therefore, the following corollaries are no longer proven here.

\newtheorem{corollary}{\hspace{2em}Corollary}[section]
\begin{corollary}
In the above theorem, if the condition $x_1<x_2$ is replaced by the condition $x_1>x_2$, it is still proved.
\end{corollary}

\newtheorem{corollary3}{\hspace{2em}Corollary}[section]
\begin{corollary}
In the above theorem, if the condition $f(x_1)<f(x_2)$ is replaced by the condition $f(x_1)>f(x_2)$, it is still proved.
\end{corollary}

\newtheorem{corollary4}{\hspace{2em}Corollary}[section]
\begin{corollary}
In the above theorem, if the condition $f(x_1)<f(x_2)$ is replaced by the condition $f(x_1)=f(x_2)$, it is still proved with 2 zero points for $T(h,x)$.
\end{corollary}

\newtheorem{corollary5}{\hspace{2em}Corollary}[section]
\begin{corollary}
In the above theorem, if two minimum values are replaced by two maximum values, it is still proved.
\end{corollary}

\newtheorem{corollary6}{\hspace{2em}Corollary}[section]
\begin{corollary}
In the above theorem, if the function is on the interval of $[A, B]$ rather than the whole define domain, it is still proved.
\end{corollary}

\newtheorem{corollary7}{\hspace{2em}Corollary}[section]
\begin{corollary}
In the above theorem, the positive and negative symbols of $T(h_0,x)$ remain unchanged within the interval $[x_0-\frac{h_0}{2},x_0+\frac{h_0}{2}]$ corresponding to that there is no extreme point $x_0$. And $T(h_0,x)$ is always greater than or equal to $0$ or less than $0$.
\end{corollary}

Next, a further theorem is proved.
The average gradient of the points around the extreme point which is not the global extreme point must not be $0$, but there is one average gradient of the points around the global extreme point must be $0$.
\newtheorem{theorem5}{\hspace{2em}Theorem}[section]
\begin{theorem}
For the function $f(x)$, the corresponding zero point of the average gradient function $T(h_0,x)$ is $x_T$.
The corresponding interval $[x_T-\frac{h_0}{2},x_T+\frac{h_0}{2}]$ with the scale $h$ contains the extreme point of $f(x)$, and the extreme point is the global extreme point, except points on the border.
\end{theorem}

\newtheorem{proof5}{\hspace{2em}Proof}
\begin{proof}
$\because T(h_0,x_T)=\frac{f(x_T+\frac{h_0}{2})-f(x_T-\frac{h_0}{2})}{h}=0$, and $h>0$.

$\therefore f(x_T+\frac{h_0}{2})-f(x_T-\frac{h_0}{2})=0$.

According to the Lagrange mean value theorem, $\exists \xi \in [x_T-\frac{h_0}{2},x_T+\frac{h_0}{2}]$, $f'(\xi)=\frac{f(x_T+\frac{h_0}{2})-f(x_T-\frac{h_0}{2})}{h_0}=0$.

$\therefore \xi$ is the extreme point of function $f(x)$.

$\because$ For the interval corresponding to the extreme point of function $f(x)$, $\forall x$, $T(h_0,x)$ has no sign-changing zero point.

$\therefore \xi$ is the extreme point of the function $f(x)$.
\end{proof}[QED]

\newtheorem{theorem6}{\hspace{2em}Theorem}[section]
\begin{theorem}
For the function $f(x)$ in the above theorem, the zero point $x_T$ of its corresponding average gradient function $T(h_0,x)$ is on the two sides of the global extreme point $x_0$ of $f(x)$, and between extreme points or between boundary values.
\end{theorem}

\newtheorem{proof6}{\hspace{2em}Proof}
\begin{proof}
Take one case of the two as an example.

$\because x_T \in [x_A,x_3]$, $f(x_A)>f(x_1)$, $f(x_3)>f(x_1)$.

$\therefore T(h_0,x_A+\frac{h_0}{2})=\frac{f(x_A+h_0)-f(x_A)}{h_0}<0$, $T(h_0,x_3-\frac{h_0}{2})=\frac{f(x_3-h_0)-f(x_3)}{h_0}>0$.

$\therefore$ The situation theorem is proved.

The other situations are similarly proved.
\end{proof}[QED]

Then, we can get the following inference.
\newtheorem{corollary8}{\hspace{2em}Corollary}[section]
\begin{corollary}
If the global extreme point of the function $f(x)$ in the interval $[A,B]$ is on the boundary, $T(h_0,x)$ is always not less than $0$ or not greater than $0$ in the interval.
\end{corollary}

Finally, give a general theorem.
\newtheorem{theorem7}{\hspace{2em}Theorem}[section]
\begin{theorem}
The function $f(x)$ is continuous and derivable on $R$, and all extreme points are not equal.
Then the average gradient function $T(h,x)$ can be constructed on the interval $[A,B]$ (which can be $[\infty,-\infty]$) by appropriately selecting the scale $h$ value.
As $h$ increases, the zero points of $T(h,x)$ becomes less and less.
The interval corresponding to each zero point of $T(h,x)$ contains at least one extreme point of the function $f(x)$.
When $h$ increases to a certain value, i.e. $h \ge h_0$ (similar to the $h_0$ value in the above theorem, but may be smaller), $T(h,x)$ will have only 2, 1 or 0 zero point(s), or is constant $0$.
Specifically (without considering maximum or minimum points at boundary),

(1) If the global maximum point and the global minimum point of $f(x)$ are extreme points, $T(h_0,x)$ will have 2 zero points, and for the zero point, $T(h_0,x+\triangle x )T(h_0,x-\triangle x) < 0$ (where $\triangle x$ is infinitely small).

(2) If $f(x)$ only has the global maximum (or minimum) point as the extreme point, then $T(h_0,x)$ will have 1 zero point, and for the zero point, $T(h_0,x+ \triangle x)T(h_0,x-\triangle x) < 0$ (where $\triangle x$ is infinitely small).

(3) If $f(x)$ is equal to a constant value, $T(h_0,x)$ is always $0$.
\end{theorem}

\newtheorem{proof7}{\hspace{2em}Certificate}
\begin{proof}
Using the above theorems to compare each adjacent two extreme points,then this theorem is proved.
\end{proof}[QED]

\newtheorem{corollary9}{\hspace{2em}inference}[section]
\begin{corollary}
In the above theorem, if the function $f(x)$ has $m$ equal maximum points and $n$ equal minimum points.
Then the corresponding conclusion becomes,

(1) If the global maximum points and the minimum points of $f(x)$ are extreme points, $T(h_0,x)$ will have $(m+n)$ zero points, and for the zero points, $T(H_0,x+\triangle x)T(h_0,x-\triangle x) < 0$ (where $\triangle x$ is infinitely small).

(2) The situation involved boundary points is similar to the theorem above.
\end{corollary}

\section{Homogenization Function}

In the proof of the theorem of the previous section, only the average gradient function of the function is used.
It is easy to think that if take the integral of this average gradient function, then it is a good estimate of the function.
And according to the above theorem, the integral of this average gradient function should estimate the monotonicity of most parts of a function.
So it also screened out most of the extreme points that are not the maximum (or minimum) points.

Before the next analysis, this paper will demonstrate from another perspective what the average gradient function does.
The average gradient function $T(h,x)$ of a unary function is equivalent to convolution of the derivative function, and the convolution kernel function is
$${
h(x)=
\begin{cases}
1 &t-\frac{h}{2}<x<t+\frac{h}{2}\\
0 &\text{Others}
\end{cases}.
}$$
$T(h,x)=\frac{1}{h}\int h(t)f(x-t)dt=\frac{1}{h}\int_{x-\frac{h}{2}} ^{x+\frac{h}{2}}f(t)dt $.
This homogenization makes the function smoother.
Integral the average gradient function $T(h,x)$ by $\int T(h,x)dx$ to get an estimate of the function $f(x)$, which is called the homogenization function.

\newtheorem{definition5}{\hspace{2em}Definition}[section]
\begin{definition}
$\int T(h,x)dx$ is called the {\bf homogenization function} of the function $f(x)$.
\end{definition}
Using the homogenization function of the unary function, it is also possible to prove the theorems of the previous section and obtain the corresponding theorems.
The purpose of extracting the homogenization function is for the next argument.
Multivariate functions must rely on homogenization functions to obtain a more concise form.
Because the multidimensional characteristics of the multivariate function make the derivative a number of partial derivatives.

\section{Average Gradient and Homogenization Function of Multivariate Function}

First, discuss the continuous and derivable binary function $f(x,y)$ as an example.
According to the definition of the unary function average gradient function, the binary function $f(x, y)$ is the average of all derivative values on a two-dimensional region $[A_x, B_x; A_y, B_y]$.
What corresponding to the derivative of the unary function are two partial derivatives, so the average gradient function is a vector $(T_x(h,x,y), T_y(h,x,y))$, where
$$
T_x(h,x,y)=T_x=\frac{1}{h^2}\int_{y-\frac{h}{2}}^{y+\frac{h}{2}} \int_{ X-\frac{h}{2}}^{x+\frac{h}{2}}f_x(t_x,t_y)dt_xdt_y,
$$
$$
T_y(h,x,y)=T_y=\frac{1}{h^2}\int_{y-\frac{h}{2}}^{y+\frac{h}{2}} \int_{ X-\frac{h}{2}}^{x+\frac{h}{2}}f_y(t_x,t_y)dt_xdt_y.
$$
Therefore, a vector field is obtained according to the definition of the average gradient, which becomes average gradient field.

\newtheorem{definition6}{\hspace{2em}Definition}[section]
\begin{definition}
The vector field, consisting of $T_x(h,x,y)$ and $T_y(h,x,y)$ defined by partial derivatives, is called the {\bf average gradient field} of function $f(x,y)$.
\end{definition}

It can be proved that the vector field is a potential field because $\frac{\partial T_x}{\partial y} = \frac{\partial T_y}{\partial x}$.

\newtheorem{definition7}{\hspace{2em}Definition}[section]
\begin{definition}
The average gradient field of the function $f(x,y)$ is a potential field, and its potential function is the {\bf homogenization function},
\begin{equation}
F(h,x,y)=\frac{1}{h^2}\int_{y-\frac{h}{2}}^{y+\frac{h}{2}} \int_{x- \frac{h}{2}}^{x+\frac{h}{2}}f(t_x,t_y)dt_xdt_y,
\end{equation}
of the function $f(x,y)$.
\end{definition}
This can be proved and found according to the method of potential function of potential field, in order to verify whether the homogenization function $F(h, x, y)$ and the average gradient field have similar properties to those of the unary function.
Prove a lemma first.

\newtheorem{lemma}{\hspace{2em} Lemma}[section]
\begin{lemma}
For the binary function $f(x,y)$, substitute $y$ with $ax+b$ then get $g(x)=f(x,ax+b)$.
For any $a$ and $b$, $g(x)$ is a gradient solvable function (a function that can be got the maximum (or minimum) points by the gradient algorithm. Note that this is not necessarily a convex function, it can be a quasi-convex function).
Then, $f(x,y)$ is a gradient solvable function.
\end{lemma}

\newtheorem{proof8}{\hspace{2em} Proof}
\begin{proof}

This actually decomposes the gradient algorithm.
Discuss the case of the gradient descent method.

First take a point $(x_1,y_1)$ randomly, and take a line $y=a_1x+b_1$, passing through $(x_1,y_1)$.

$\because$ $g(x)=f(x,a_1x+b_1)$ is a gradient solvable function.

$\therefore$ The minimum point $x_2$ is obtained by $g(x)$, and $a_1x+b_1$ is substituted to get $y_2=a_1x_2+ b_1$, and get the point $(x_2,y_2)$.

Then take a different line $y=a_2x+b_2$ passing through $(x_2,y_2)$.

Repeating the above steps, it is able to find the minimum point of $f(x,y)$.

$\because$ $(x_i,y_i)$ is arbitrary, and each new point found is based on the gradient descending.

$\therefore$ $f(x,y)$ is the gradient solvable function.

The situation for the gradient ascending method is similar.
\end{proof}[QED]

\newtheorem{corollary11}{\hspace{2em}inference}[section]
\begin{corollary}
The minimum (or maximum) point of the binary function $f(x,y)$ in the above theorem is $(x_0, y_0)$.
Substitute $y$ with $ax+b$ to get $g(x)=f(x,ax+b)$.
If $\forall a,b$, the minimum (or maximum) point of $g(x)$ is $(x_g,ax_g+b)$, which is $(x_g,y_g)$.
Then the point $(x_g,y_g)$ is the point $(x_0, y_0)$.
The anti-proposition is also established.
\end{corollary}

If $f(x,y)$ is understood as a surface, then the above lemma says that an arbitrary curve intercepted by a line is gradient solvable, then $f(x,y)$ is gradient solvable.
The dot column $\{(x_i,y_i)\}$ is convergent and converges to the global extreme point.

\newtheorem{theorem9}{\hspace{2em}Theorem}[section]
\begin{theorem}
$$
F(h,x,y)=\frac{1}{h^2}\int_{y-\frac{h}{2}}^{y+\frac{h}{2}} \int_{x- \frac{h}{2}}^{x+\frac{h}{2}}f(t_x,t_y)dt_xdt_y
$$
is the homogenization function of the binary function $f(x,y)$. The function $f(x,y)$ is continuous and derivable.
There is a positive value of $h_0$ such that $F(h_0,x,y)$ is a gradient solvable function, and the minimum value of $F(h_0,x,y)$ is $(x_T,y_T)$ corresponding to the region $[x_T-\frac{h_0}{2},x_T+\frac{h_0}{ 2};y_T-\frac{h_0}{2},y_T+\frac{h_0}{2}]$ contains the global extreme point $(x_0,y_0)$ of $f(x,y)$.
\end{theorem}

\newtheorem{proof9}{\hspace{2em}Certificate}
\begin{proof}
According to the theorem of the unary function, there is a positive value of $h_1$, and $y=a_1x+b_1$ is substituted into $F(h,x,y)$ to get $G(h_1,x)$.
$G(h_1,x)$ is a gradient solvable function.
$$
G(h_1,x)=F(h_1,x,a_1x+b_1)=\frac{1}{h_1^2}\int_{a_1x+b_1-\frac{h_1}{2}}^{a_1x+b_1+\frac{h_1}{2}} \int_{x-\frac{h_1}{2}}^{x+\frac{h_1}{2}}f(t_x,a_1t_x+b_1)dt_xd(a_1t_x+b_1).
$$
There may be different $h$ values for different lines $y=ax+b$.

$\because a_1$, $b_1$ is arbitrary, according to the lemma, $F(h_1,x,y)$ is a gradient solvable function.
And the minimum $h_0$ for $F(h,x,y)$ is the maximum value of all the $h_i$ values that may be obtained in the above operation.

$\therefore$ There is one $h_0$ value, which makes $F(h,x,y)$ a gradient solvable function.

$\therefore F(h,x,y)$ has a minimum (or maximum) point $(x_T,y_T)$, satisfying $\frac{\partial F}{\partial x}(h_0,x_T,y_T)=0$ and $\frac{\partial F}{\partial y}(h_0,x_T,y_T)=0$.

$\therefore \frac{1}{h_0^2}\int_{y_T-\frac{h_0}{2}}^{y_T+\frac{h_0}{2}}\int_{x_T-\frac{h_0}{ 2}}^{x_T+\frac{h_0}{2}}f_x(t_x,t_y)dt_xdt_y=0$,
$\frac{1}{h_0^2}\int_{y_T-\frac{h_0}{2}}^{y_T+\frac{h_0}{2}}\int_{x_T-\frac{h_0}{2} }^{x_T+\frac{h_0}{2}}f_y(t_x,t_y)dt_xdt_y=0$.

Let the area $D=[x_T-\frac{h}{2}, x_T+\frac{h}{2}; y_T-\frac{h}{2}, y_T+\frac{h}{2}]$.

$\therefore \iint f_x(t_x,t_y)dt_xdt_y=0$,
$\iint f_y(t_x,t_y)dt_xdt_y=0$.

$\because$ The function $f(x,y)$ is continuous and derivable in the area $D$.

$\therefore$ There is at least one point $(\xi_{x_1}, \xi_{y_1}) \in D$, and $f_x(\xi_{x_1},\xi_{y_1})=0$.

For the same reason, there is at least one point $(\xi_{x_2}, \xi_{y_2}) \in D$, and $f_y(\xi_{x_2},\xi_{y_2})=0$.

$\because f(x,y)$ is continuous in the area $D$.

$\therefore$ There is at least one point $(\xi_x,\xi_y) \in D$, and $f_x(\xi_x, \xi_y)=f_y(\xi_x, \xi_y)=0$.

$\therefore$ There has the extreme point $(\xi_x,\xi_y) \in D$ of $f(x,y)$.

$\because$ Based on the previous $F(h_0,x,y)$ search method, for any $y=ax+b$ passing through $(\xi_{x_0}, \xi_{y_0})$, substitute $y$ with $ax+b$ in $f(x,y)$ and $F(h_0,x,y)$ to get $g(x)$ and $G(h_0,x)$.

$\therefore$ According to the conclusion of the unary function, the interval corresponding to the minimum point of $G(h,x)$ is $[x_T-\frac{h}{2},x_T+\frac{h}{2}]$ Contains the minimum point of $g(x)$.

$\because$ For any $a$ and $b$, the curve $g(x)$ obtained from $f(x,y)$ has a minimum point of $( \xi_x,\xi_y)$.

$\therefore$ According to the lemma inference, $(\xi_x,\xi_y)$ is the maximum value of the function $f(x,y)$.
\end{proof}[QED]

\newtheorem{theorem10}{\hspace{2em}Theorem}[section]
\begin{theorem}

The function $f(x,y)$ is continuous and derivable on $R$ and all function values of extreme points are different.
By choosing a appropriate positive value of $h$ in the region $[A_x, B_x; A_y, B_y]$, an homogenization function $F(h, x, y)$ can be constructed.
With the increase of $h$, $F(h,x,y)$ has fewer and fewer extreme points.
The regions of $f(x,y)$ corresponding to the extreme points of $F(h,x,y)$ have extreme points of $f(x,y)$.
When $h$ increases to a certain value, i.e. $h>h_0$ (similar to the $h_0$ in the above theorem, but possibly smaller).
$F(h,x,y)$ will have only 2, 1 or 0 extreme point(s) or is constant 0.
Specifically (without considering maximum or minimum points at boundary),

(1) If both the global maximum point and the global minimum point of $f(x,y)$ are extreme points, then $F(h_0,x,y)$ will have two extreme points, which are the global maximum point and the global minimum point respectively.

(2) If only the global maximum point (or the global minimum point) of $f(x, y)$ is the extreme point, then $F(h_0, x, y)$ will have an extreme point, which is the global maximum (or minimum) point.

(3) If the maximum point of $f(x,y)$ is at the interval point, $F(h_0,x,y)$ has no extreme point or is constant.

(4) If $f(x, y)$ equals a constant value, $F(h_0, x, y)$ is constant.
\end{theorem}

\newtheorem{corollary12}{\hspace{2em}Corollary}[section]
\begin{corollary}

In the above theorem, if the function $f(x,y)$ has $m(m \ge 1)$ equal global maximum points and $n(n \ge 1)$ equal global minimums points.
Then the corresponding conclusion becomes:

(1) If the global maximum point(s) and global minimum point(s) of $f(x,y)$ are extremum points, then $F(h_0,x,y)$ will have m global maximum points and n global minimum points without other extreme point.

(2) If only global maximum point(s) (or global minimum point(s)) of $f(x, y)$ are the extreme point, then $F(h_0,x,y)$ will have m(or n) global maximum points (or global minimum points) without other extreme point.

\end{corollary}

\newtheorem{corollary13}{\hspace{2em}Corollary}[section]
\begin{corollary}
In the above theorem, if the value of $h$ is very large, then the function $F(h,x,y)$ will become a plane, a saddle surface or a surface similar to saddle surface in a finite region.
\end{corollary}

The case of multivariate functions can be obtained similarly according to the above theorem.
The homogenization function of the multivariate function $f(x_1,x_2,\cdots,x_n)$ is
\begin{equation}
F(h,x_1,x_2,\cdots,x_n)=\frac{1}{h^n}\int_{x_n-\frac{h}{2}}^{x_n+\frac{h}{2}} \cdots \int_{x_1-\frac{h}{2}}^{x_1+\frac{h}{2}} f(x_1,x_2,\cdots,x_n)dx_1\cdots dx_n.
\end{equation}

All these proofs imply a method for changing all functions into convex optimization solvable functions.

\end{document}